\documentclass{amsart}
\usepackage{amsthm,amsmath,amssymb,amsfonts,enumerate,amscd}
\DeclareMathAlphabet{\mathscr}{OT1}{pzc}%
                                 {m}{it}
\DeclareMathOperator\rad{rad}
\DeclareMathOperator\soc{soc}

\DeclareMathOperator\Hom{Hom}
\DeclareMathOperator\End{End}
\DeclareMathOperator\Ext{Ext}
\DeclareMathOperator\Mod{Mod}

\def\rMod#1{\mathscr{Mod}_{\hspace{-2pt} #1\hspace{-0pt}}}

\def\FinLen#1{\mathscr{FL}_{\hspace{-1pt} #1\hspace{-0pt}}}

\newtheorem{theorem}{Theorem}
\newtheorem{proposition}[theorem]{Proposition}

\newtheorem{lemma}[theorem]{Lemma}
\newtheorem{question}[theorem]{Question}

\theoremstyle{definition}
\newtheorem{definition}[theorem]{Definition}
\newtheorem{example}[theorem]{Example}

\theoremstyle{remark}

\newtheorem{case'}{Case}

\numberwithin{theorem}{section}

\begin{document}
\title{Extensions of Simple Modules and the Converse of Schur's Lemma}
\author{Greg Marks}
\address{Department of Mathematics and Computer Science\\
St.\ Louis University\\
St.\ Louis, MO 63103--2007\\
U.S.A.}
\email{marks@slu.edu}
\author{Markus Schmidmeier}
\address{Mathematical Sciences\\
Florida Atlantic University\\
Boca Raton, FL 33431--0991\\
U.S.A.}
\email{markus@math.fau.edu}
\keywords{Converse of Schur's Lemma, Gabriel quiver}
\subjclass{Primary 16D90, 16G20, 16S50}
\date{May 13, 2009}

\begin{abstract}
The {\em converse of Schur's lemma} (or {\em CSL}) condition on a module
category has been the subject of considerable study in recent years.  In
this note we extend that work by developing basic properties of module
categories in which the CSL condition governs modules of finite length.
\end{abstract}

\maketitle

\section{Introduction}\label{Introduction}

Schur's Lemma states that for any ring $R$ and any simple module $M_R$,
the endomorphism ring $\End(M_R)$ is a division ring.  In this note we
are interested in the converse of Schur's Lemma (CSL), i.e.\ whether for
a given module category $\mathscr{C},$ every object in $\mathscr{C}$
whose endomorphism ring is a division ring is in fact simple.  If this
is the case, we say that $\mathscr{C}$ {\em has CSL}.  The case that has
received almost exclusive attention in the literature (e.g.\ see
\cite{ah1}, \cite{ah2}, \cite{f}, \cite{hk}, \cite{hp}, \cite{wz}) is
$\mathscr{C} = \rMod{R}$, the category of right $R$-modules.  Here we
will focus on the case $\mathscr{C} = \FinLen{R}$, the category of right
$R$-modules of finite length.

We propose to separate the study of rings $R$ which satisfy CSL for
$\FinLen{R}$ from the study of rings which satisfy CSL for $\FinLen{R}$
but not for $\Mod_R$, since the two properties relate to different
topics: extensions of simples versus constructions of large modules.

It turns out that the CSL property for finite length modules---and
sometimes the CSL property for all modules---is controlled by the
following combinatorial information:

\begin{definition}

Let $R$ be a ring.  Recall that the {\em right Gabriel quiver} (or {\em
right $\Ext$-quiver}) of $R$ is the directed graph $Q$ consisting of the
following data:

 \begin{itemize}

 \item The points of $Q$ are in bijective correspondence with the
 isomorphism classes of simple right $R$-modules.

 \item There is an arrow $i \rightarrow j$ in $Q$ whenever the 
 corresponding simple modules $S_i$ and $S_j$ extend, i.e.\ 
 $\Ext_R^1(S_i,S_j) \neq \{0\}$.

 \end{itemize}

\noindent
We say that the right Gabriel quiver is {\em totally disconnected} if
there are no arrows between any two different points.

\end{definition}

For example, the right Gabriel quiver of a semisimple ring is a disjoint
union of finitely many points.  The right Gabriel quiver of the discrete
valuation ring $R = k[x]_{(x)}$ (for $k$ a field) has one point and one
loop.  

The reader should be aware that the literature contains some variants of
the definition we give here.  In the classical setting where $R$ is a
finite-dimensional algebra over a field, some authors adopt the
convention that in the right Gabriel quiver of $R$ the arrow between the
vertices corresponding to $S_i$ and $S_j$ carries as label the pair
given by the dimensions of $\Ext_R^1(S_i,S_j)$ as a vector space over
$\End(S_j)_R$ and $\End(S_i)_R$ respectively.

\begin{theorem}\label{thm-fl}
Let $R$ be any ring.  Then $\FinLen{R}$ has CSL if and only if
the right Gabriel quiver of $R$ is totally disconnected.
\end{theorem}

We will prove Theorem~\ref{thm-fl} in Section~\ref{sect-fl}.

In general, for $\FinLen{R}$ to have CSL is a considerably weaker
condition than for $\rMod{R}$ to have CSL\@. The distinction between CSL
on $\FinLen{R}$ and CSL on $\rMod{R}$ is illustrated in the following
examples.

\begin{example}\label{commutative_finite_length_csl}
Let $R$ be any commutative ring whatsoever.  Then $\FinLen{R}$ has
CSL\@.  To see this, suppose $M$ is an $R$-module of finite length such
that $\End(M_R)$ is a division ring.  If $M$ were not simple, then by
\cite[Corollary]{wz}, $M_R$ would be isomorphic to the field of
fractions of $R/\mathfrak{p}$ where $\mathfrak{p} =
\operatorname{ann}^R(M)$ is a prime but not maximal ideal of $R$,
contradicting the hypothesis that $M$ has finite length.

By contrast, in \cite{wz} it is shown that for a commutative ring $R$,
the category $\rMod{R}$ has CSL if and only if $R$ has Krull dimension
$0$.
\end{example}

We infer from Theorem~\ref{thm-fl} that the (right) Gabriel quiver of
any commutative ring is totally disconnected.

Example~\ref{commutative_finite_length_csl} suggests a further reason
why the (not necessarily commutative) rings $R$ for which $\FinLen{R}$
has CSL are an interesting object of study: they include all commutative
rings, so this condition is a new sort of generalization of
commutativity.

\begin{example}[J. H. Cozzens \cite{coz}]\label{ex-coz}
Let $K$ be an algebraically closed field of positive characteristic
$p$, let $\varphi$\mbox{\rm :}~$x\mapsto x^{p^n}$ be a Frobenius
automorphism on $K$, and let $k = K^{\langle\varphi\rangle}$ be the
fixed field of $\varphi$. Then the skew Laurent polynomial ring
$R=K[x,x^{-1};\varphi]$ has, up to isomorphism, a unique simple right
module $S$, which is injective. Thus, every finite length right $R$-module
is semisimple, and hence $\FinLen{R}$ has CSL\@. 

Nevertheless, $\rMod{R}$ does not have CSL\@.  It is easy to show
that if $R$ is a right or left Ore domain, then $\rMod{R}$ has CSL if
and only if $R$ is a division ring.  In the present example $R$ is a
simple noetherian domain, hence an Ore domain.

Note that the category $\FinLen{R}$ is equivalent to the category
$\FinLen{k}$, which obviously has CSL (as even $\rMod{k}$ does).
\end{example}

We can generalize this example using \cite[Theorem A]{schm}.  Let
$\rad(A)$ denote the Jacobson radical of a ring $A$.  Recall that a
finite-dimensional $k$-algebra $A$ is called {\em elementary} if
$A/\!\rad(A)$ is a finite direct product of copies of $k$.

\begin{proposition}\label{prop-elementary}
Let $A$ be a finite-dimensional elementary algebra over a finite field
whose right Gabriel quiver is totally disconnected.  Then there exists a
noetherian ring $R$ such that the following properties hold:
\begin{enumerate}
\item[\rm(i)] The categories $\FinLen{A}$ and $\FinLen{R}$ are
equivalent.
\item[\rm(ii)] Both $\FinLen{A}$ and $\FinLen{R}$ have CSL\@.
\item[\rm(iii)] The category $\rMod{A}$ has CSL, but the category
$\rMod{R}$ does not.
\end{enumerate}
\end{proposition}

On the other hand, for semiprimary rings the CSL property for all
modules is controlled by the right Gabriel quiver, i.e.\ it is
controlled by the CSL property for finite length modules.  Recall that a
ring $R$ is said to be {\em semiprimary} if the Jacobson radical
$\rad(R)$ is nilpotent and $R/\!\rad(R)$ is a semisimple ring.
Semiprimary rings figure prominently in our main object of study here:
it is well known that the endomorphism ring of a finite length module is
semiprimary.

\begin{theorem}\label{thm-semiprimary}
Let $R$ be a semiprimary ring.  The following conditions are
equivalent:
\begin{enumerate}
\item[\rm(i)] $\FinLen{R}$ has CSL.
\item[\rm(ii)] The right Gabriel quiver of $R$ is totally disconnected.
\item[\rm(iii)] $R$ is a finite direct product of full matrix rings
over local rings.
\item[\rm(iv)] $\rMod{R}$ has CSL\@.
\end{enumerate}
\end{theorem}

We defer the proofs of Proposition~\ref{prop-elementary} and
Theorem~\ref{thm-semiprimary} to Section~\ref{sect-all}.  The literature
contains results akin to Theorem~\ref{thm-semiprimary}, such as the
following.

\begin{theorem}\label{thm-noetherian}
Let $R$ be a one-sided noetherian ring or a perfect ring.  Then
$\rMod{R}$ has CSL if and only if $R$ is a finite direct product
of full matrix rings over local perfect rings.
\end{theorem}

The left noetherian case is covered by \cite[Theorem~1]{ah2}, the right
noetherian case by \cite[Theorem~3.4]{dm}, and the perfect case by
\cite[Theorem~1.2]{ah1}.

\begin{example}\label{ex-weyl}
Let $k$ be a field of characteristic $0$, and let $\mathbb{A}_1\!(k)
= k\langle x, y \rangle / (xy - yx - 1)$ be the first Weyl algebra
over $k$.  If $$S_1 = \mathbb{A}_1\!(k) / x \mathbb{A}_1\!(k)
\qquad\mbox{and}\qquad S_2 = \mathbb{A}_1\!(k) / (x+y) \mathbb{A}_1\!(k),$$
then by \cite[Proposition~5.6, Theorem~5.7]{mr}, $S_1$ and $S_2$
are nonisomorphic simple right $\mathbb{A}_1\!(k)$-modules for which
$\Ext_{\mathbb{A}_1\!(k)}^1(S_1, S_2) \neq \{0\}$.  Therefore, the right
Gabriel quiver of $\mathbb{A}_1\!(k)$ is not totally disconnected,
so Theorem~\ref{thm-fl} tells us $\FinLen{\mathbb{A}_1\!(k)}$ does
not have CSL\@.
\end{example}

The conclusion of Example~\ref{ex-weyl} can be extended from
$\mathbb{A}_1\!(k)$ to certain {\em generalized Weyl algebras;} see
\cite[Theorem~1.1]{bavula} for details.

\begin{example}\label{ex-bounded_dpr}
Let $R$ be a right bounded Dedekind prime ring.  Then $\FinLen{R}$
has CSL\@.  To prove this, first note that if $R$ is right primitive
then by \cite[Theorem~4.10]{er} it is simple artinian (so in this
case even $\rMod{R}$ has CSL)\@.  Now assume $R$ is not right
primitive.  Suppose $S_1$ and $S_2$ are arbitrary nonisomorphic
simple right $R$-modules.  Then $\operatorname{ann}^R_r(S_1) =
\mathfrak{m}_1$ and $\operatorname{ann}^R_r(S_2) = \mathfrak{m}_2$
are maximal ideals of $R$ and $\mathfrak{m}_1 \neq \mathfrak{m}_2$,
by \cite[Theorem~3.5]{robson}.  By \cite[Theorem~1.2, Proposition~2.2]{er},
$\mathfrak{m}_1$ and $\mathfrak{m}_2$ are invertible ideals.  We
can therefore apply \cite[Proposition~1]{gw} to conclude that
$\Ext_{\mathbb{A}_1\!(k)}^1(S_1, S_2) = \{0\}$.  Thus, by
Theorem~\ref{thm-fl}, $\FinLen{R}$ has CSL\@.
\end{example}

\begin{example}\label{ex-group_algebras}
Let $G$ be a finite group and $F$ a field.
\begin{enumerate}

\item[\rm(i)]
If the characteristic of $F$ does not divide the order of $G$, then
$FG$ is semisimple and hence the Gabriel quiver is totally disconnected
(no arrows).

\smallskip

\item[\rm(ii)]
If the characteristic of $F$ is a prime number $p$ and $G$ is a
finite $p$-group, then $FG$ is a local ring and again, the Gabriel
quiver is totally disconnected (the only arrow is a loop).

\smallskip

\item[\rm(iii)]
In the case where the number of simples is different from the number
of blocks, there is a block where two nonisomorphic simples extend,
and we get a proper arrow in the Gabriel quiver.
\end{enumerate}
Thus, in cases (i) and (ii), but not (iii), $\FinLen{FG}$ has CSL\@.
\end{example}

\section{CSL for finite length modules}\label{sect-fl}

In this section we give a proof of Theorem~\ref{thm-fl}.  First, assume
that $\FinLen{R}$ has CSL\@.  As a consequence of the next lemma, the
right Gabriel quiver of $R$ must be totally disconnected.

\begin{lemma} \label{division_subring}
Suppose $0 \rightarrow T \rightarrow M \rightarrow S \rightarrow
0$ is a non-split short exact sequence in $\rMod{R}$ where $S$ and
$T$ are nonisomorphic simple modules.  Then $\End(M_R)$ is isomorphic
to a division subring of $\End(S_R)$ and of $\End(T_R)$.
\end{lemma}

\begin{proof}
Since $\Hom_R(M,T)=\{0\}$, $\Hom_R(M,M)$ embeds in $\Hom_R(M,S)$;
since $\Hom_R(T,S)=\{0\}$, we can identify $\Hom_R(M,S)$ with
$\Hom_R(S,S)$.  This yields a ring monomorphism $\End(M_R) \rightarrow
\End(S_R)$.  Similarly, since $\Hom_R(S,M)=\{0\}$ and $\Hom_R(T,S)=\{0\}$,
we obtain a ring monomorphism $\End(M_R) \rightarrow \End(T_R)$.
Thus $\End(M_R)$ is isomorphic to a subring of the division rings
$\End(S_R)$ and $\End(T_R)$, so $\End(M_R)$ is a domain.  Being
also semiprimary, $\End(M_R)$ is a division ring.
\end{proof}

For the converse, we assume that the right Gabriel quiver of $R$ is
totally disconnected.  We first show that every finite length
indecomposable module is isotypic, and then that every isotypic module
is either simple or admits a nonzero nilpotent endomorphism.  Note that
some of the results apply both to finite length modules over an
arbitrary ring and to arbitrary modules over a semiprimary ring.  We
will use these results again in the next section.

\begin{definition}
A module $M$ is {\em isotypic} if all simple subquotients of $M$
are isomorphic.  A sequence
$$0=M_0\subset M_1\subset \cdots\subset M_\ell=M$$
of submodules is called an {\em isotypic filtration} of $M$ of
length $\ell$ if for every $i = 1, \ldots, \ell$, the quotient
$M_i/M_{i-1}$ is isotypic.
\end{definition}

\begin{proposition}\label{prop-isotypic}
Suppose $R$ is a ring whose right Gabriel quiver is totally
disconnected.  Suppose that either

\begin{enumerate}

\item[\rm(i)] $M$ is an object of $\FinLen{R}$, or

\item[\rm(ii)] $R$ is semiprimary, and $M$ is an object of $\rMod{R}$.

\end{enumerate}
Then $M$ is a finite direct sum of isotypic modules.
\end{proposition}

\begin{proof}
{\it Step 1: The module $M$ has an isotypic filtration.}  For example,
the radical filtration of $M$ can be refined to an isotypic
filtration.

{\it Step 2: For each isotypic filtration 
$$0\subset M_1\subset\cdots\subset M_\ell=M$$
there is an isotypic filtration $0\subset M_1'\subset\cdots\subset
M_\ell'=M$ such that $M_i\subseteq M_i'$ for each $i$ and $\Hom_R(M_1',
M/M_1') = \{0\}$.} Zorn's Lemma can be applied to the set of all
isotypic submodules of $M$ that contain $M_1$; let $M_1'$ be a
maximal member of this set.  We have $\Hom_R(M_1', M/M_1')=\{0\}$
since the socle of $M/M_1'$ cannot contain a simple summand isomorphic
to a subquotient of $M_1'$.  For $i>1$, put $M_i'=M_i+M_1'$. Then
$$\frac{M_i'}{M_{i-1}'} = \frac{M_i+M_1'}{M_{i-1}+M_1'} \cong
\frac{M_i}{M_i\cap(M_{i-1}+M_1')} = \frac{M_i}{M_{i-1}+(M_i\cap
M_1')}$$ is epimorphic image of $M_i/M_{i-1}$ and hence isotypic.

{\it Step 3: Let $M$ and $N$ be isotypic modules that both satisfy {\rm
(i)} or {\rm (ii)} of the proposition and for which $\Hom_R(M,N)=\{0\}$.
Then $\Ext^1_R(M,N)=\{0\}$.}  When $M$ and $N$ are semisimple, this
follows from the hypothesis on the right Gabriel quiver.  The general
case follows by induction on the lengths of semisimple filtrations of
$M$ and $N$.

{\it Step 4: The module $M$ is a direct sum of isotypic modules.}
We induct on the length $\ell$ of the isotypic filtration of $M$
produced in Step~2.  The case $\ell=1$ is trivial.  For the induction
step, let $M$ have an isotypic filtration $0\subset
M_1'\subset\cdots\subset M_{\ell+1}'=M$.  By inductive hypothesis,
$M/M_1' \cong \bigoplus_i M_i''$ is a direct sum of isotypic modules
$M_i''$.  Since $\Hom_R(M_1', M_i'') = \{0\}$ for all $i$, by Step~3
we have $$\Ext^1_R(M_1', M/M_1') = \bigoplus_i \Ext^1_R(M_1', M_i'')
= \{0\}$$ Thus, the short exact sequence $0 \rightarrow M_1'
\rightarrow M \rightarrow M/M_1' \rightarrow 0$ splits, and $M \cong
M_1\oplus\bigoplus_iM_i''$ is a direct sum of isotypic modules.
\end{proof}

Next we prove a criterion for isotypic modules to admit a nonzero
nilpotent endomorphism.  The argument is adapted from the proof of
\cite[Theorem~1.2]{ah1}.

\begin{lemma}\label{lemma-soc}
Let $R$ be any ring, and let $M$ be a right $R$-module.  Then $M$
has no nonzero semisimple direct summand if and only if $\soc(M)
\subseteq \rad(M)$.
\end{lemma}
\begin{proof}
If $M$ has no nonzero semisimple direct summand, then every simple
submodule is superfluous, whence $\soc(M) \subseteq \rad(M)$.
Conversely, if $M$ does have a nonzero semisimple direct summand,
then $M$ has a simple direct summand, which is contained in $\soc(M)$
but not $\rad(M)$, so $\soc(M) \not\subseteq \rad(M)$.
\end{proof}

\begin{proposition}\label{prop-nilpotent}
Suppose that $M_R$ is a nonzero isotypic module that is not simple.
Assume in addition that either $M$ has finite length or $R$ is a
perfect ring.  Then $M$ has a nonzero nilpotent endomorphism.
\end{proposition}

\begin{proof}
If $M$ has a nonzero simple direct summand, the conclusion is clear;
so assume otherwise.  By Lemma~\ref{lemma-soc}, $\soc(M) \subseteq
\rad(M)$.  Now, $M$ is nonzero and isotypic, and the hypotheses
imply that $M/\!\rad(M)$ is semisimple; therefore, there exists a
nonzero homomorphism $f_0$\mbox{\rm :}~$M/\!\rad(M) \rightarrow
\soc(M)$.  The composite map $$\mbox{$f$: } M
\stackrel{\pi}{\longrightarrow} M/\!\rad(M) \stackrel{f_0}{\longrightarrow}
\soc(M) \stackrel{\iota}{\longrightarrow} M$$ (where $\pi$ is the
canonical epimorphism and $\iota$ the inclusion map) is a nonzero
endomorphism of $M$ satisfying $f^2=0$.
\end{proof}

Theorem~\ref{thm-fl} is now established.  The ``only if'' part follows
from Lemma~\ref{division_subring}.  The ``if'' part follows from
Propositions~\ref{prop-isotypic} and \ref{prop-nilpotent}.

\section{CSL for all modules}\label{sect-all}

To prove Theorem~\ref{thm-semiprimary} we will show $$\mbox{(i)}
\Leftrightarrow \mbox{(ii)} \Rightarrow \mbox{(iii)} \Rightarrow
\mbox{(iv)} \Rightarrow \mbox{(i)}.$$ By Theorem~\ref{thm-fl},
statements (i) and (ii) are equivalent for any ring $R$.

\smallskip 

$\mbox{(ii)} \Rightarrow \mbox{(iii)}$: According to
Proposition~\ref{prop-isotypic}, the module $R_R=\bigoplus_iP_i$ is a
finite direct sum of indecomposable isotypic submodules $P_i$.  Two such
submodules are either isomorphic or have no nonzero homomorphisms
between them.  Thus, $R=\End(R_R)$ is a finite direct product of matrix
rings over the local endomorphism rings of the $P_i$'s.

\smallskip 

$\mbox{(iii)} \Rightarrow \mbox{(iv)}$: Apply ``$\mbox{(iv)} \Rightarrow
\mbox{(iii)}$'' of \cite[Theorem~1.2]{ah1}.
\qed

\smallskip

\medskip

We now prove Proposition~\ref{prop-elementary}.  Let $A$ be an
elementary algebra over a finite field $k$ of $p^n$ elements.  Let $K$
be an algebraically closed field of characteristic $p$ and
$\varphi$\mbox{\rm :}~$K\to K$ the Frobenius automorphism, given by $\alpha
\mapsto \alpha^{p^n}$; we identify $k$ with the fixed field of
$\varphi$.  Let $\Sigma=K[x,x^{-1};\varphi]$ be the V-ring studied in
\cite{coz}; we claim that the ring $R=\Sigma\otimes_kA$ has the required
properties.

\smallskip 

(i) We infer from \cite[Theorem A]{schm} that the categories
of $\FinLen{A}$ and $\FinLen{R}$ are equivalent.

\smallskip

(ii) Follows from (i) and (iii).

\smallskip

(iii) Applying Theorem~\ref{thm-semiprimary} to the semiprimary
ring $A$, we deduce that $\rMod{A}$ has CSL\@.  To see that $R$
does not have CSL, note that since $A$ is elementary, the ring
homomorphism $\pi$\mbox{\rm :}~$A \rightarrow k$ gives rise to a
surjective ring homomorphism $$\mbox{$\pi\otimes 1$: } \,
R=A\otimes_k\Sigma \longrightarrow k\otimes_k\Sigma.$$ Since
$\rMod{\Sigma}$ does not have CSL (as explained in Example~\ref{ex-coz}),
and $\Sigma$ is isomorphic to a factor ring of $R$, $\rMod{R}$ does
not have CSL\@.
\qed

\medskip

\section{Some questions}\label{sect-r2}

The rings in Examples~\ref{ex-coz} and \ref{ex-weyl} are both simple
noetherian domains.  In light of the diametrically different behavior in
the two examples, we pose the following question.

\begin{question} 
For which simple noetherian domains $R$ does $\FinLen{R}$ have CSL?
\end{question} 

\begin{question}
When do other subcategories of $\rMod{R}$ have CSL\@?  What are the
conditions under which all artinian modules have CSL\@?   All noetherian
modules?  Is there an example of a ring which has CSL for finite
length modules, but not CSL for artinian modules?
\end{question}

One may also consider categories with quasi-CSL in the following sense:

\begin{definition}
Let $\mathscr{C}$ be a category of modules, i.e.\ $\mathscr{C}$ is
a full subcategory of $\rMod{R}$ for some ring $R$.  We say that
an object $M$ of $\mathscr{C}$ is {\em quasi-simple} if the only
submodules $N\subseteq M$ such that $N$ and $M/N$ are objects of
$\mathscr{C}$ are $N=0$ and $N=M$.  The category $\mathscr{C}$ is
said to have {\em quasi-CSL} if the only modules with endomorphism
ring a division ring are the quasi-simple ones.
\end{definition}

\begin{question}
Are there interesting categories with quasi-CSL?
\end{question}

\begin{example}
For a given ring $R$, the category $\FinLen{R}$ or the category
$\rMod{R}$ has quasi-CSL if and only if it has CSL\@.
\end{example}

Nevertheless, in general {\it quasi-CSL} and {\it CSL} are different
conditions on a module category, as will be seen in
Example~\ref{quasi-CSL-ex} below.  We preface this example with some
motivating observations.

\begin{example}
Even if $R$ has CSL, then the ring $U_2(R)$ of upper triangular $2$
by $2$ matrices with coefficients in $R$ need not have CSL\@.
Indeed, if $S$ is a simple right $R$-module, then the row
$\left(\begin{matrix} S & S \end{matrix}\right)$ is a right module
over $U_2(R)$ of length $2$ with endomorphism ring $\End(S_R)$.

In fact, the category $\rMod{U_2(R)}$ is just the category of all
maps between right $R$-modules, a map $f$\mbox{\rm :}~$A \rightarrow
B$ being a module over $U_2(R)$ via $$(a,b)\cdot {\textstyle
\left(\begin{matrix} x & y \\ 0 & z \end{matrix}
\right)}\;=\;(ax,f(a)y+bz).$$ Consider the full subcategory $\mathcal
S(R)$ of $\rMod{U_2(R)}$ consisting of all maps which are monomorphisms
(``$\mathcal S$=submodules'').
\end{example}

\begin{question}
If $R$ has CSL, does $\mathcal S(R)$ have quasi-CSL?
\end{question}

Categories of type $\mathcal S(R)$ play a role in applications of
ring theory; for example, the embeddings of a subgroup in a
finite abelian group, or the embeddings of a subspace in a vector
spaces such that the subspace is invariant under the action of a
linear operator, fall into this type of category.

\begin{example} \label{quasi-CSL-ex}
For $\Lambda$ a commutative uniserial ring with radical generator
$p$ and radical factor field $k$, the category $\mathcal S(\Lambda)$
has quasi-CSL but not CSL, as follows.  There are exactly two
quasi-simple modules, $S_1 = \left(\begin{matrix} k & k
\end{matrix}\right)$ and $S_2 = \left(\begin{matrix} 0 & k
\end{matrix}\right)$, up to isomorphy; both have endomorphism ring
$k$.  Then $\mathcal S(\Lambda)$ has quasi-CSL since any embedding
$\left(\begin{matrix} A & B \end{matrix}\right)$ with $B$ a semisimple
$\Lambda$-module is a direct sum of copies of $S_1$ and $S_2$.  On
the other hand, if $B$ is not semisimple then multiplication
by $p$ is a nonzero nilpotent endomorphism.
\end{example}

\def\auth#1{{\rm #1,}}
\def\titlart#1{{\it #1,}}
\def\titlj#1{{\rm #1\mbox{}}}
\def\titlbook#1{{\it #1\mbox{}}}
\def\vol#1{{\bf #1}}
\def\no#1{{\rm no.\ #1,}}
\def\date#1{{\rm (#1),}}
\def\pages#1{{\rm #1.}}


\begin{thebibliography}{99}

\bibitem{ah1} 
\auth{M. Alaoui, A. Haily}
\titlart{Perfect rings for which the converse of Schur's lemma
holds}
\titlj{Publ.\ Mat.} 
\vol{45}
\date{2001}
\no{1}
\pages{219--222}

\bibitem{ah2} 
\auth{M. Alaoui, A. Haily}
\titlart{The converse of Schur's lemma in noetherian rings and group
algebras}
\titlj{Comm.\ Algebra}
\vol{33}
\date{2005}
\no{7}
\pages{2109--2114}

\bibitem{bavula}
\auth{V. V. Bavula}
\titlart{The extension group of the simple modules over the first Weyl
algebra}
\titlj{Bull.\ London Math.\ Soc.}
\vol{32}
\date{2000}
\no{2}
\pages{182--190}

\bibitem{coz} 
\auth{J. H. Cozzens}
\titlart{Homological properties of the ring of differential
polynomials}
\titlj{Bull.\ Amer.\ Math.\ Soc.}
\vol{76}
\date{1970}
\pages{75--79}

\bibitem{dm}
\auth{M. Dombrovskaya, G. Marks}
\titlart{Asymmetry in the converse of Schur's Lemma}
\titlj{Comm.\ Algebra,}
to appear.

\bibitem{er}
\auth{D. Eisenbud, J. C. Robson}
\titlart{Hereditary Noetherian prime rings}
\titlj{J. Algebra}
\vol{16}
\date{1970}
\pages{86--104}

\bibitem{f}
\auth{C. Faith}
\titlart{Indecomposable injective modules and a theorem of Kaplansky}
\titlj{Comm.\ Algebra}
\vol{30}
\date{2002}
\no{12}
\pages{5875--5889}

\bibitem{gw}
\auth{K. R. Goodearl, R. B. Warfield Jr.}
\titlart{Simple modules over hereditary Noetherian prime rings}
\titlj{J. Algebra}
\vol{57}
\date{1979}
\no{1}
\pages{82--100}

\bibitem{hp} 
\auth{Y. Hirano, J. K. Park}
\titlart{Rings for which the converse of Schur's Lemma holds}
\titlj{Math.\ J. Okayama Univ.}
\vol{33}
\date{1991}
\pages{121--131}

\bibitem{hk}
\auth{C. Huh, C. O. Kim} 
\titlart{$\pi$-regular rings satisfying the converse of Schur's lemma}
\titlj{Math.\ J. Okayama Univ.}
\vol{34}
\date{1992}
\pages{153--156}

\bibitem{mr}
\auth{J. C. McConnell, J. C. Robson}
\titlart{Homomorphisms and extensions of modules over certain
differential polynomial rings}
\titlj{J. Algebra}
\vol{26}
\date{1973}
\pages{319--342}

\bibitem{robson}
\auth{J. C. Robson}
\titlart{Non-commutative Dedekind rings}
\titlj{J. Algebra}
\vol{9}
\date{1968}
\pages{249--265}

\bibitem{schm} 
\auth{M. Schmidmeier}
\titlart{A family of noetherian rings with their finite length
modules under control}
\titlj{Czechoslovak Math.\ J.}
\vol{52(127)}
\date{2002}
\no{3}
\pages{545--552}

\bibitem{wz}
\auth{R. Ware, J. Zelmanowitz}
\titlart{Simple endomorphism rings}
\titlj{Amer.\ Math.\ Monthly}
\vol{77}
\date{1970}
\no{9}
\pages{987--989}

\end{thebibliography}
\end{document}